\documentclass[11pt,a4paper]{article}
\usepackage[T1]{fontenc}
\usepackage[utf8]{inputenc}
\usepackage{lmodern}
\usepackage[a4paper,margin=27mm]{geometry}
\usepackage{amsmath,amssymb,amsthm,mathtools}
\usepackage{booktabs,array}
\usepackage{graphicx}
\usepackage{algorithm,algpseudocode}
\usepackage{enumitem}
\usepackage{microtype}
\usepackage[hidelinks]{hyperref}
\hypersetup{pdftitle={Approximate Sampling of Doubly Free-Boundary Macdonald Processes},pdfsubject={Reflection, fixed-size resampling, and total variation convergence},pdfauthor={}}
\numberwithin{equation}{section}
\newtheorem{theorem}{Theorem}[section]
\newtheorem{proposition}[theorem]{Proposition}
\newtheorem{lemma}[theorem]{Lemma}
\newtheorem{corollary}[theorem]{Corollary}
\theoremstyle{definition}
\newtheorem{definition}[theorem]{Definition}

\newcommand{\Y}{\mathbb Y}
\newcommand{\N}{\mathbb Z_{\geq 0}}
\newcommand{\PP}{\mathbb P}

\newcommand{\ind}{\mathbf 1}
\newcommand{\el}{\mathrm{el}}
\newcommand{\oa}{\mathrm{oa}}
\newcommand{\deel}{\mathrm{deel}}
\newcommand{\eoa}{\mathrm{eoa}}
\newcommand{\TV}{\mathrm{TV}}
\newcommand{\Geom}{\operatorname{Geom}}
\newcommand{\Bern}{\operatorname{Bern}}
\newcommand{\NB}{\mathsf{NB}}
\newcommand{\ZZ}{\mathcal Z}
\newcommand{\seedZ}{\mathfrak Z}
\newcommand{\cset}{\mathcal C}
\newcommand{\size}[1]{\lvert #1\rvert}
\newcommand{\tv}[1]{\left\lVert #1\right\rVert_{\TV}}

\setlist[enumerate]{itemsep=3pt,topsep=5pt}
\allowdisplaybreaks[2]
\title{\textbf{A Randomized Approximation Algorithm for\\
Doubly Free-Boundary Macdonald Processes}}
\author{
  Zhongyang Li\\[4pt]
  {\small Department of Mathematics}\\
  {\small University of Connecticut}\\
  {\small\url{https://mathzhongyangli.wordpress.com}}
}
\date{}

\begin{document}
\maketitle
\vspace{-1.6em}
\begin{abstract}
We construct an approximate sampler for finite doubly free-boundary
Macdonald rail-yard processes with Macdonald, boundary, and edge
parameters in $(0,1)$. The construction allows arbitrary finite
sequences of additions and removals of horizontal or vertical strips
and all four boundary types. It reverses a finite
commutation--reflection history from a partition sequence whose
entries are equal, with the common partition weighted by both
boundaries. We derive local reconstruction kernels from the skew
Cauchy, Littlewood, and branching identities. Bulk and boundary
updates combine integer increments with conditional shape sampling
at a prescribed size; same-sign swaps use finite conditional resampling.

For $K$ reflection periods and an initial partition-size cutoff $M$,
the algorithm terminates almost surely and samples the resulting
approximating law exactly. We prove total variation convergence to
the target process as $K,M\to\infty$ and give explicit bounds for both
truncation errors, decaying geometrically in the cutoffs for fixed
input parameters. Under fixed parameter bounds, we also bound the
expected operation count of direct enumeration in a unit-cost model.

We further construct an exact rejection sampler for the initial
partition law without a size cutoff or evaluation of its normalizing
constant. Using this sampler removes the size-truncation error,
leaving only the error from the finite reflection history.
\end{abstract}

\section{Introduction and scope}

We construct an explicit approximate sampler for the finite doubly
free-boundary Macdonald rail-yard processes of Li and Shi
\cite{LiShi2026}. The sampler applies to general Macdonald parameters,
arbitrary finite sequences of additions and removals of horizontal
or vertical strips, and all four boundary types. It generates the
entire partition sequence, with explicit total variation bounds
on its distance from the target law.

The algorithm starts from a constant partition sequence and reconstructs
the original path by reversing a finite sequence of commutations and
boundary reflections. The common starting partition is random:
even when all edge parameters vanish, it need not be empty.
Its law retains the weights of both boundaries and must therefore
be sampled with those weights.

The local reconstruction also retains the full Macdonald shape weights.
In a bulk or boundary update, the current partitions and a sampled
integer increment determine the size of the new partition. We then
sample its shape from the allowed partitions of that size, using
their Macdonald weights. A same-sign swap instead samples the middle
partition from the finite interval between its fixed outer endpoints.
These kernels provide stochastic reconstruction compatible with
both free boundaries.

For $K$ reflection periods and an initial partition-size cutoff $M$,
the algorithm terminates almost surely and samples an explicitly
defined approximating law exactly. We prove convergence to the target
process in total variation as $K,M\to\infty$. We bound the reflection
error and the initial size-truncation error separately; for each fixed
input, both bounds decay geometrically in their respective cutoffs.

We also construct an exact rejection sampler for the weighted starting
partition without a size cutoff or evaluation of its normalizing
constant. Using this sampler removes the size-truncation error,
while the finite reflection history still produces an approximation.
For the implementation with a size cutoff, we give an expected
operation-count bound for direct enumeration under fixed parameter
bounds in a unit-cost model.

Several ingredients have precedents. Fulman
\cite[Section~5]{Fulman1997} samples ordinary Macdonald measures by
first drawing the size of a horizontal strip and then drawing its
shape with normalized Pieri weights. Betea et al.\
\cite{SchurSample2018} construct exact Schur-process samplers from
bijective Cauchy and Littlewood identities and growth diagrams.
Barraquand, Borodin, and Corwin
\cite[Sections~2.4.1--2.4.2]{BBC2020} construct bulk and boundary
Markov kernels for growing half-space Macdonald processes at general
$(q,t)$. Their push-block kernels sample shapes using normalized
Macdonald and Littlewood weights. Our construction uses this
weight-transport principle in a history involving repeated reflections
at two free boundaries, together with the resulting weighted starting
law and explicit truncation bounds.

Li's Schur sampler \cite[Section~6, Theorem~6.9]{Li2023}
uses repeated commutations and reflections, samples finite truncations
exactly, and proves total variation convergence. Li and Shi
\cite[Section~3, Proposition~3.14]{LiShi2026} establish the corresponding
Macdonald identities and reduction scheme. The present paper turns
that reduction into a complete sampling construction with the full
Macdonald weights, implements the weighted starting law, and gives
quantitative error and enumeration-cost bounds. Even in the Schur
specialization $q=t$, our sampler is in general different from that
of \cite{Li2023}: its bulk and boundary reconstruction uses conditional
shape resampling at a prescribed size rather than the bijective local
rules employed there; see Section~\ref{sec:discussion}.

We consider a process with a fixed number $N\geq1$ of transitions.
Here $q,t$ are the Macdonald parameters, $u,v$ are the left and right
boundary parameters, respectively, and $x_i$ is the edge parameter of
the $i$-th transition. We assume
\begin{equation}\label{eq:regime}
0<q,t,u,v<1,\qquad 0<x_i<1\quad(1\leq i\leq N).
\end{equation}
These conditions ensure that the distributions used by the sampler
are well defined and that the target partition function is finite.
A proof of finiteness is given in Appendix~\ref{app:growth}.

\section{The process and its scalar factors}\label{sec:model}

\subsection{Partitions, branching coefficients, and dual letters}

Let $\Y$ be the set of integer partitions, including $\varnothing$. We identify partitions with their Young diagrams.
For $\mu\subseteq\lambda$, the skew diagram $\lambda/\mu$
consists of the boxes of $\lambda$ that do not belong to $\mu$.

Write $\size{\lambda}=\sum_i\lambda_i$, and let $\lambda'$ be the conjugate partition. For $s=(i,j)\in\lambda$, set
\[
a_\lambda(s)=\lambda_i-j,\qquad
\ell_\lambda(s)=\lambda'_j-i,
\qquad
b_\lambda(s)=\frac{1-q^{a_\lambda(s)}t^{\ell_\lambda(s)+1}}
{1-q^{a_\lambda(s)+1}t^{\ell_\lambda(s)}}.
\]
Define $b_\lambda=\prod_{s\in\lambda}b_\lambda(s)$, with $b_\varnothing=1$. 

Let $P_\lambda=P_\lambda(\,\cdot\,;q,t)$ denote the monic Macdonald
symmetric function, and set $Q_\lambda=b_\lambda P_\lambda$.
We write $P_{\lambda/\mu}$ and $Q_{\lambda/\mu}$ for the corresponding
skew Macdonald symmetric functions, using the conventions of
\cite[Chapter~VI]{Macdonald1995}. They satisfy
\begin{equation}\label{eq:skew-ratio}
Q_{\lambda/\mu}
=\frac{b_\lambda}{b_\mu}P_{\lambda/\mu}.
\end{equation}

For a type $a\in\{L,R\}$, introduce the nonnegative one-parameter specializations
\begin{equation}\label{eq:specializations}
p_k(\rho_L(x))=x^k,
\qquad
p_k(\rho_R(x))=(-1)^{k-1}\frac{1-q^k}{1-t^k}x^k.
\end{equation}
The term ``nonnegative'' refers to their skew Macdonald evaluations, not to the signs of their power sums. Define
\begin{align}
U_a(\lambda/\mu;x)&:=P_{\lambda/\mu}(\rho_a(x);q,t),\nonumber\\
V_a(\lambda/\mu;x)&:=Q_{\lambda/\mu}(\rho_a(x);q,t).\label{eq:UV}
\end{align}
Equivalently, using the duality for skew Macdonald functions
\cite[Chapter~VI, (7.16)]{Macdonald1995}, we have
\begin{equation}\label{eq:UV-table}
\begin{array}{c|cc}
 & U_a(\lambda/\mu;x)&V_a(\lambda/\mu;x)\\ \hline
L&P_{\lambda/\mu}(x;q,t)&Q_{\lambda/\mu}(x;q,t)\\[2pt]
R&Q_{\lambda'/\mu'}(x;t,q)&P_{\lambda'/\mu'}(x;t,q).
\end{array}
\end{equation}
For $\mu\subseteq\lambda$, the skew diagram $\lambda/\mu$ is called
a \emph{horizontal strip} if it contains at most one box in each
column, and a \emph{vertical strip} if it contains at most one box
in each row. The weights $U_L(\lambda/\mu;x)$ and
$V_L(\lambda/\mu;x)$ vanish unless $\lambda/\mu$ is a horizontal
strip, while $U_R(\lambda/\mu;x)$ and $V_R(\lambda/\mu;x)$ vanish
unless $\lambda/\mu$ is a vertical strip. Homogeneity gives coefficients $u_a,v_a$ such that
\begin{equation}\label{eq:uv-coeff}
U_a(\lambda/\mu;x)=u_a(\lambda/\mu)x^{\size{\lambda}-\size{\mu}},\qquad
V_a(\lambda/\mu;x)=v_a(\lambda/\mu)x^{\size{\lambda}-\size{\mu}}.
\end{equation}
Coefficients are zero unless the indicated strip is allowed, and positive otherwise. In particular,
\begin{equation}\label{eq:coeff-ratio}
v_a(\lambda/\mu)=\frac{b_\lambda}{b_\mu}u_a(\lambda/\mu).
\end{equation}
Appendix~\ref{app:branching} gives finite product formulas for all these coefficients. The symmetric-function conventions and identities used here are standard; see \cite[Chapter~VI]{Macdonald1995}.

\subsection{Boundary weights and the target measure}

Set
\[
b^\el_\lambda
=\prod_{\substack{s\in\lambda\\
                  \ell_\lambda(s)\text{ even}}}b_\lambda(s),
\qquad
b^\oa_\lambda
=\prod_{\substack{s\in\lambda\\
                  a_\lambda(s)\text{ odd}}}b_\lambda(s).
\]
and let $b^{\deel}_\lambda=b^\el_\lambda$, $b^{\eoa}_\lambda=b^\oa_\lambda$. For $c\in\cset:=\{\el,\oa,\deel,\eoa\}$ define
\[
\ind_c(\lambda)=
\begin{cases}
1,&c=\el,\oa,\\
\ind_{\{\lambda'\text{ has only even parts}\}},&c=\deel,\\
\ind_{\{\lambda\text{ has only even parts}\}},&c=\eoa,
\end{cases}
\]
and abbreviate
\begin{equation}\label{eq:beta-tau}
\beta_c(\lambda)=\ind_c(\lambda)b^c_\lambda,
\qquad
\tau_c(\lambda)=\frac{\beta_c(\lambda)}{b_\lambda}.
\end{equation}
The left boundary uses $\tau_c$, and the right boundary uses $\beta_c$.

A labelled word is a sequence
\[
w=((a_1,\sigma_1,x_1),\ldots,(a_N,\sigma_N,x_N)),
\qquad a_i\in\{L,R\},\quad\sigma_i\in\{+,-\}.
\]

We regard the $i$-th entry $(a_i,\sigma_i,x_i)$, together with
its fixed label $i$, as a labelled letter. A labelled letter is
called a minus letter or a plus letter according as its sign
$\sigma_i$ is $-$ or $+$. Adjacent swaps move the entire
labelled letter, including its type, sign, and parameter.

For $\alpha,\gamma\in\Y$, define the transition weights
\[
E_{a,+}(\alpha,\gamma;x)=U_a(\gamma/\alpha;x),\qquad
E_{a,-}(\alpha,\gamma;x)=V_a(\alpha/\gamma;x).
\]
For $\boldsymbol\lambda=(\lambda^{(0)},\ldots,\lambda^{(N)})\in\Y^{N+1}$ and boundary types $c_l,c_r\in\cset$, define
\begin{equation}\label{eq:target-weight}
W_w(\boldsymbol\lambda)=
\tau_{c_l}(\lambda^{(0)})u^{\size{\lambda^{(0)}}}
\beta_{c_r}(\lambda^{(N)})v^{\size{\lambda^{(N)}}}
\prod_{i=1}^N E_{a_i,\sigma_i}
 (\lambda^{(i-1)},\lambda^{(i)};x_i).
\end{equation}
The target process is
\begin{equation}\label{eq:target-law}
\pi_w(\boldsymbol\lambda)=\frac{W_w(\boldsymbol\lambda)}{Z_w},
\qquad Z_w=\sum_{\boldsymbol\lambda\in\Y^{N+1}}W_w(\boldsymbol\lambda).
\end{equation}
Zero weights enforce all interlacing and boundary restrictions. 

\subsection{Bulk factors}

For $0<Q<1$ and $0\leq A<1$, use the notation
\begin{equation}\label{eq:k-def}
K_{A,Q}(r)=\frac{(Ar;Q)_\infty}{(r;Q)_\infty},\qquad
k_n(A,Q)=\frac{(A;Q)_n}{(Q;Q)_n},
\qquad K_{A,Q}(r)=\sum_{n\geq0}k_n(A,Q)r^n.
\end{equation}
Here $(A;Q)_n=\prod_{j=0}^{n-1}(1-AQ^j)$ and $k_0=1$. The final equality is the $q$-binomial identity \cite[\S17.5]{DLMF}; it also follows by comparing coefficients in
$(1-r)K_{A,Q}(r)=(1-Ar)K_{A,Q}(Qr)$.

Consider two adjacent steps $(a,+,x)$ and $(b,-,y)$,
where $a,b\in\{L,R\}$. Keep the outer partitions fixed and sum
over the intermediate partition. By the skew Cauchy identity,
the total weight in this order is $C_{ab}(xy)$ times the total
weight in the reversed order $(b,-,y),(a,+,x)$.
Writing $r=xy$, we have
\begin{equation}\label{eq:C-factors}
C_{ab}(r)=\sum_{n\geq0}c_n^{ab}r^n=
\begin{cases}
K_{t,q}(r),&a=b=L,\\
K_{q,t}(r),&a=b=R,\\
1+r,&a\ne b.
\end{cases}
\end{equation}
Consequently, the $RR$ factor uses the \emph{swapped} parameters. This follows directly from \eqref{eq:specializations} and the Cauchy exponential
\[
\Pi(\rho,\sigma;q,t)=\exp\left(\sum_{k\geq1}
\frac{1-t^k}{1-q^k}\frac{p_k(\rho)p_k(\sigma)}{k}\right).
\]
We will always compute the factors from these explicit specialization conventions.

\subsection{Boundary factors and their coefficients}\label{subsec:F}

Let
\[
\Theta_c(X;q,t)=\sum_{\lambda\in\Y}\beta_c(\lambda)P_\lambda(X;q,t),
\qquad F_{c,a}(r)=\Theta_c(\rho_a(r);q,t).
\]
The Littlewood identities underlying the reflections are recorded in
\cite[Lemma~3.2 and Corollaries~3.6, 3.8]{LiShi2026}. Their single-variable factors and coefficients are
\begin{equation}\label{eq:F-table}
\begin{array}{c|c|c}
c & F_{c,L}(r) & d_n^{c,L}:=[r^n]F_{c,L}(r)\\ \hline
\el
& K_{t,q}(r)
& k_n(t,q)\\[3pt]
\oa
& (1+r)K_{qt,q^2}(r^2)
& k_{\lfloor n/2\rfloor}(qt,q^2)\\[3pt]
\deel
& 1
& \ind_{\{n=0\}}\\[3pt]
\eoa
& K_{qt,q^2}(r^2)
& \begin{cases}
  k_{n/2}(qt,q^2),&n\text{ even},\\
  0,&n\text{ odd}.
  \end{cases}
\end{array}
\end{equation}
In particular,
\[
F_{\oa,L}(r)=\frac{(qtr^2;q^2)_\infty}
{(1-r)(q^2r^2;q^2)_\infty},
\]
so its coefficients of $r^{2m}$ and $r^{2m+1}$ agree.

Define the involution on boundary types by
\begin{equation}\label{eq:hat}
\widehat\el=\oa,\qquad\widehat\oa=\el,
\qquad\widehat\deel=\eoa,\qquad\widehat\eoa=\deel.
\end{equation}
Then
\begin{equation}\label{eq:F-dual}
F_{c,R}(r;q,t)=F_{\widehat c,L}(r;t,q),
\qquad d_n^{c,R}(q,t)=d_n^{\widehat c,L}(t,q).
\end{equation}
Indeed, conjugation interchanges arms and legs and gives
\[
b_{\lambda'}(t,q)=b_\lambda(q,t)^{-1},\qquad
\beta_{\widehat c}(\lambda';t,q)
=\frac{\beta_c(\lambda;q,t)}{b_\lambda(q,t)}.
\]
Substitute these relations and \eqref{eq:UV-table} into the series for $\Theta_c$ to obtain \eqref{eq:F-dual}. 

\section{Scalar sampling without infinite products}\label{sec:scalar}

For later use define the distribution
\begin{equation}\label{eq:NB}
\NB_Q(A,r)\{n\}=\frac{k_n(A,Q)r^n}{K_{A,Q}(r)},
\qquad n\in\N,\quad 0<r<1.
\end{equation}
This notation specifies the parameter convention completely. Bulk increments are independent of the current partitions and have laws
\begin{equation}\label{eq:bulk-increments}
G_{ab}(r)\sim
\begin{cases}
\NB_q(t,r),&a=b=L,\\
\NB_t(q,r),&a=b=R,\\
\Bern\bigl(r/(1+r)\bigr),&a\ne b.
\end{cases}
\end{equation}
A boundary reflection contributes the scalar factor
$F_{c,a}(r)=\sum_{n\geq0}d_n^{c,a}r^n$.
We sample the associated integer increment by assigning probability
proportional to $d_n^{c,a}r^n$ to the value $n$.
Thus, for $0<r<1$, define the nonnegative-integer random variable
$B_{c,a}(r)$ by
\begin{equation}\label{eq:boundary-increments}
\PP(B_{c,a}(r)=n)
=\frac{d_n^{c,a}r^n}{F_{c,a}(r)},
\qquad n\in\mathbb Z_{\geq0}.
\end{equation}

For $a=L$, the following constructions implement this law:
\begin{align}
B_{\el,L}(r)&\sim\NB_q(t,r),\nonumber\\
B_{\oa,L}(r)&=2J+\epsilon,
&&J\sim\NB_{q^2}(qt,r^2),\quad
\epsilon\sim\Bern\bigl(r/(1+r)\bigr),\label{eq:boundary-samplers}\\
B_{\deel,L}(r)&=0,\nonumber\\
B_{\eoa,L}(r)&=2J,
&&J\sim\NB_{q^2}(qt,r^2).\nonumber
\end{align}
In the second line $J$ and $\epsilon$ are independent. For $a=R$, applying \eqref{eq:F-dual} gives
\begin{align*}
B_{\el,R}(r)&\overset{\mathrm d}{=}2J+\epsilon,
&&J\sim\NB_{t^2}(qt,r^2),\quad
  \epsilon\sim\Bern\bigl(r/(1+r)\bigr),\\
B_{\oa,R}(r)&\sim\NB_t(q,r),\\
B_{\deel,R}(r)&\overset{\mathrm d}{=}2J,
&&J\sim\NB_{t^2}(qt,r^2),\\
B_{\eoa,R}(r)&=0.
\end{align*}
In the first line, $J$ and $\epsilon$ are independent.

Although \eqref{eq:NB} contains an infinite-product normalizer, the following rejection sampler does not need it. Our convention is
$\PP(\Geom(r)=n)=(1-r)r^n$.

\begin{lemma}[A finite-product scalar sampler]\label{lem:rejection}
Let $L_Q=\lceil(1-Q)^{-1}\rceil$ and $H_Q=(1-Q)^{-L_Q}$. Repeatedly draw $J\sim\Geom(r)$ and accept $J$ with probability $k_J(A,Q)/H_Q$. The accepted value has law \eqref{eq:NB}; the procedure terminates almost surely. Each acceptance probability uses only a finite product.
\end{lemma}
\begin{proof}
Since $(A;Q)_n\leq1$,
\[
k_n(A,Q)\leq\frac1{(Q;Q)_\infty}.
\]
Moreover,
\begin{align*}
\log\frac1{(Q;Q)_\infty}
&=\sum_{m\geq1}\frac{Q^m}{m(1-Q^m)}
\leq\frac{-\log(1-Q)}{1-Q}
\leq L_Q[-\log(1-Q)].
\end{align*}
Thus $k_n\leq H_Q$. The probability of proposing and accepting $n$ is
$(1-r)r^n k_n/H_Q$, which normalizes to \eqref{eq:NB}. The acceptance probability per trial is
$(1-r)K_{A,Q}(r)/H_Q>0$, proving almost-sure termination.
\end{proof}

\section{Exact local shape kernels}\label{sec:kernels}

Given the current partitions, each scalar increment determines the
size of the updated partition. We now specify the conditional
distribution on the allowed shapes of that size and prove that
the resulting local updates reproduce the required Macdonald weights.

Coefficients with negative indices are understood to be zero.

\subsection{Bulk peak reconstruction}

Fix the outer partitions $\lambda,\mu$, a $+$ type $a$ with parameter $x$, and a $-$ type $b$ with parameter $y$. The skew Cauchy identity is
\begin{equation}\label{eq:skew-cauchy}
\sum_\nu U_a(\nu/\lambda;x)V_b(\nu/\mu;y)
=C_{ab}(xy)\sum_\kappa V_b(\lambda/\kappa;y)U_a(\mu/\kappa;x).
\end{equation}
It changes the local word $(a,+,x),(b,-,y)$ into $(b,-,y),(a,+,x)$.
The left side contains a peak $\nu$ and the right side a valley $\kappa$.
The identity follows from the usual skew Cauchy identity after applying
\eqref{eq:specializations}; compare \cite[Lemma~3.9]{LiShi2026}.

Define, for $m\in\N$,
\begin{equation}\label{eq:Sbulk}
S_m^{a,b}(\lambda,\mu)=\sum_{\size{\zeta}=m}
 u_a(\zeta/\lambda)v_b(\zeta/\mu).
\end{equation}
Comparing homogeneous coefficients in \eqref{eq:skew-cauchy} gives
\begin{equation}\label{eq:bulk-coeff-id}
S_m^{a,b}(\lambda,\mu)=
\sum_\kappa v_b(\lambda/\kappa)u_a(\mu/\kappa)
c_{m+\size{\kappa}-\size{\lambda}-\size{\mu}}^{ab}.
\end{equation}
Here $c_n^{ab}$ are the coefficients defined in
\eqref{eq:C-factors}, with $c_n^{ab}=0$ for $n<0$. The right sum is finite: a nonzero term has $\kappa\subseteq\lambda\cap\mu$.

\begin{definition}[Bulk reconstruction kernel]\label{def:bulk-kernel}
Given an admissible valley $\kappa$, draw $G\sim G_{ab}(xy)$ independently and set
\begin{equation}\label{eq:bulk-size}
m=\size{\lambda}+\size{\mu}-\size{\kappa}+G.
\end{equation}
Then sample a partition $\nu$ according to
\begin{equation}\label{eq:bulk-shape}
\PP(\nu\mid\lambda,\mu,\kappa,G)=
\frac{\ind_{\{\size{\nu}=m\}}u_a(\nu/\lambda)v_b(\nu/\mu)}
{S_m^{a,b}(\lambda,\mu)}.
\end{equation}
Let $T_{a,b}^{x,y}(\kappa,\nu\mid\lambda,\mu)$ denote the
probability that the above procedure produces $\nu$ from $\kappa$,
with the outer partitions $\lambda,\mu$ fixed. For an admissible
input $\kappa$, put
\[
g=|\nu|+|\kappa|-|\lambda|-|\mu|.
\]
Then
\[
T_{a,b}^{x,y}(\kappa,\nu\mid\lambda,\mu)
=
\begin{cases}
\displaystyle
\frac{c_g^{ab}(xy)^g}{C_{ab}(xy)}
\frac{u_a(\nu/\lambda)v_b(\nu/\mu)}
     {S_{|\nu|}^{a,b}(\lambda,\mu)},
& g\geq0 \text{ and } c_g^{ab}>0,\\[8pt]
0, & \text{otherwise}.
\end{cases}
\]
\end{definition}

The normalizing sum $S_m^{a,b}(\lambda,\mu)$ in
\eqref{eq:Sbulk} has finitely many terms, since there are only
finitely many partitions of size $m$. If the input coefficient and the sampled coefficient $c_G^{ab}$ are positive, then \eqref{eq:bulk-coeff-id} makes its denominator strictly positive.

\begin{proposition}[Bulk weight transport]\label{prop:bulk-transport}
The kernel in Definition~\ref{def:bulk-kernel} satisfies
\begin{align}
&\sum_\kappa V_b(\lambda/\kappa;y)U_a(\mu/\kappa;x)
 T_{a,b}^{x,y}(\kappa,\nu\mid\lambda,\mu)\nonumber\\
&\hspace{35mm}=\frac{U_a(\nu/\lambda;x)V_b(\nu/\mu;y)}{C_{ab}(xy)}.
\label{eq:bulk-balance}
\end{align}
\end{proposition}
\begin{proof}
Fix $\nu$ and put $m=\size{\nu}$. The sums over $\kappa$ below are restricted to partitions satisfying
$v_b(\lambda/\kappa)u_a(\mu/\kappa)>0$.

If $u_a(\nu/\lambda)v_b(\nu/\mu)=0$, both sides of
\eqref{eq:bulk-balance} vanish by the definition of $T_{a,b}^{x,y}$.
Otherwise, $S_m^{a,b}(\lambda,\mu)>0$ by \eqref{eq:Sbulk}.

For such an input $\kappa$, put
\[
g_\kappa=m+\size{\kappa}-\size{\lambda}-\size{\mu}.
\]
The corresponding monomial satisfies
\[
y^{\size{\lambda}-\size{\kappa}}
x^{\size{\mu}-\size{\kappa}}(xy)^{g_\kappa}
=
x^{m-\size{\lambda}}y^{m-\size{\mu}}.
\]
Using the explicit formula for $T_{a,b}^{x,y}$ and homogeneity,
we obtain
\begin{align*}
&\sum_\kappa
V_b(\lambda/\kappa;y)U_a(\mu/\kappa;x)
T_{a,b}^{x,y}(\kappa,\nu\mid\lambda,\mu)\\
&\quad=
\frac{U_a(\nu/\lambda;x)V_b(\nu/\mu;y)}
     {C_{ab}(xy)S_m^{a,b}(\lambda,\mu)}
\sum_\kappa
v_b(\lambda/\kappa)u_a(\mu/\kappa)c_{g_\kappa}^{ab}\\
&\quad=
\frac{U_a(\nu/\lambda;x)V_b(\nu/\mu;y)}
     {C_{ab}(xy)}.
\end{align*}
Here coefficients with negative indices are zero, and the last
equality follows from \eqref{eq:bulk-coeff-id}.
\end{proof}

\subsection{Boundary reconstruction}

Fix $c\in\cset$. Let $X=(x_1,\ldots,x_d)$ and
$Y=(y_1,\ldots,y_e)$ be finite lists of variables, and write
$X\cup Y$ for their concatenation. We suppress the parameters
$q,t$ in $\Theta_c$.

The Littlewood product identities
\cite[Lemma~3.2]{LiShi2026} give
\begin{equation}\label{eq:theta-product}
\Theta_c(x_1,\ldots,x_d)
=
\prod_{i=1}^d F_{c,L}(x_i)
\prod_{1\leq i<j\leq d}K_{t,q}(x_ix_j).
\end{equation}
Since
\[
\Pi(X,Y;q,t)
=
\prod_{i=1}^d\prod_{j=1}^e K_{t,q}(x_i y_j),
\]
separating the factors involving variables from $X$, from $Y$,
and from both gives
\begin{equation}\label{eq:theta-coproduct}
\Theta_c(X\cup Y)
=
\Theta_c(X)\Theta_c(Y)\Pi(X,Y;q,t).
\end{equation}

Fix $\eta\in\Y$. The skew Littlewood identities
\cite[Lemmas~3.3--3.5]{LiShi2026} take the form
\[
\sum_\nu\beta_c(\nu)P_{\nu/\eta}(Y)
=
\Theta_c(Y)
\sum_\kappa\beta_c(\kappa)Q_{\eta/\kappa}(Y).
\]
They follow by comparing the coefficients of $P_\eta(X)$ in
\eqref{eq:theta-coproduct}, using the branching and skew Cauchy
identities.

For the right-boundary reflection, $\eta$ is the partition
adjacent to the boundary. Substitute $Y=\rho_a(vx)$ and
multiply both sides by $v^{\size{\eta}}$.
Using $\Theta_c(\rho_a(vx))=F_{c,a}(vx)$ and the homogeneity
relations
\begin{align*}
v^{\size{\eta}}U_a(\nu/\eta;vx)
&=v^{\size{\nu}}U_a(\nu/\eta;x),\\
v^{\size{\eta}}V_a(\eta/\kappa;vx)
&=v^{\size{\kappa}}V_a(\eta/\kappa;v^2x),
\end{align*}
we obtain
\begin{equation}\label{eq:right-reflection}
\sum_\nu\beta_c(\nu)v^{\size{\nu}}U_a(\nu/\eta;x)
=
F_{c,a}(vx)
\sum_\kappa\beta_c(\kappa)v^{\size{\kappa}}
V_a(\eta/\kappa;v^2x).
\end{equation}
Thus the last letter changes from $+$ to $-$ and $x$ changes to $v^2x$.
Using \eqref{eq:coeff-ratio} gives the left-boundary version
\begin{equation}\label{eq:left-reflection}
\sum_\nu\tau_c(\nu)u^{\size{\nu}}V_a(\nu/\eta;x)
=F_{c,a}(ux)
\sum_\kappa\tau_c(\kappa)u^{\size{\kappa}}
 U_a(\eta/\kappa;u^2x).
\end{equation}
At the left boundary the first letter changes from $-$ to $+$ and $x$ changes to $u^2x$. 

For fixed $\eta\in\Y$, define the boundary normalizing sum
\begin{equation}\label{eq:Sboundary}
S_m^{c,a}(\eta)
=
\sum_{\size{\zeta}=m}\beta_c(\zeta)u_a(\zeta/\eta).
\end{equation}

The coefficient form of \eqref{eq:right-reflection} is
\begin{equation}\label{eq:boundary-coeff-id}
S_m^{c,a}(\eta)
=
\sum_\kappa
\beta_c(\kappa)v_a(\eta/\kappa)
d_{m+\size{\kappa}-2\size{\eta}}^{c,a}.
\end{equation}

\begin{definition}[Boundary reconstruction kernel]\label{def:boundary-kernel}
For a right reflection, let $h=v$ and $c=c_r$; for a left reflection, let $h=u$ and $c=c_l$. Given the reflected boundary partition $\kappa$, its neighboring interior partition $\eta$, and the parameter $x$ \emph{before} reflection, draw
\[
B\sim B_{c,a}(hx),\qquad m=2\size{\eta}-\size{\kappa}+B.
\]
Replace the boundary partition by a sample from
\begin{equation}\label{eq:boundary-shape}
\PP(\nu\mid\eta,\kappa,B)=
\frac{\ind_{\{\size{\nu}=m\}}\beta_c(\nu)u_a(\nu/\eta)}
{S_m^{c,a}(\eta)}.
\end{equation}
\end{definition}

The denominator is positive on every positive-weight input and sampled increment, by \eqref{eq:boundary-coeff-id}. Boundary indicators in $\beta_c$ automatically enforce the even or conjugate-even restrictions.

\begin{proposition}
\label{prop:boundary-transport}
Fix $\eta\in\Y$. For the right boundary, suppose that $\kappa$
is sampled with probabilities proportional to the weights in
the sum on the right-hand side of \eqref{eq:right-reflection}.
Then the procedure in Definition~\ref{def:boundary-kernel}
produces $\nu$ with probabilities proportional to the weights
in the sum on the left-hand side.

The same statement holds for the left boundary, using
\eqref{eq:left-reflection}.
\end{proposition}
\begin{proof}
For the right boundary, fix $m=\size{\nu}$ and write
$g=m+\size{\kappa}-2\size{\eta}$. By homogeneity, the weight of $\kappa$ in the sum on the
right-hand side of \eqref{eq:right-reflection} is
\[
\beta_c(\kappa)v_a(\eta/\kappa)
v^{\size{\kappa}}(v^2x)^{\size{\eta}-\size{\kappa}}.
\]
The probability of sampling $B=g$ includes the factor $(vx)^g$.
Since $g=m+\size{\kappa}-2\size{\eta}$, we have
\[
v^{\size{\kappa}}
(v^2x)^{\size{\eta}-\size{\kappa}}
(vx)^g
=
v^m x^{m-\size{\eta}}.
\]
This is the monomial of the unreflected weight. Equation~\eqref{eq:boundary-coeff-id} then cancels the shape denominator, exactly as in Proposition~\ref{prop:bulk-transport}.

For the left boundary, the desired shape coefficient is
\[
\tau_c(\nu)v_a(\nu/\eta)
=\frac{\beta_c(\nu)}{b_\eta}u_a(\nu/\eta).
\]
The factor $b_\eta^{-1}$ is constant over possible $\nu$, so its normalized shape distribution is again \eqref{eq:boundary-shape}. On the reflected side,
\[
\tau_c(\kappa)u_a(\eta/\kappa)
=\frac{\beta_c(\kappa)}{b_\eta}v_a(\eta/\kappa).
\]
The same coefficient identity and the same monomial calculation with $u$ instead of $v$ establish the left-boundary balance.
\end{proof}

\subsection{Same-sign commutations}

For two adjacent steps with the same sign, we exchange their order,
keep the outer partitions fixed, and sample a new intermediate
partition with probabilities proportional to the product of the
two transition weights in the new order.

First consider the swap
\[
(a,+,x),(b,+,y)
\longrightarrow
(b,+,y),(a,+,x).
\]
For fixed endpoints $\alpha\subseteq\gamma$, the skew branching
identity gives
\begin{align}
&\sum_{\alpha\subseteq\nu\subseteq\gamma}
U_a(\nu/\alpha;x)U_b(\gamma/\nu;y)\nonumber\\
&\qquad=
\sum_{\alpha\subseteq\zeta\subseteq\gamma}
U_b(\zeta/\alpha;y)U_a(\gamma/\zeta;x).
\label{eq:same-sign-up}
\end{align}
Let $\zeta$ be the current intermediate partition in the swapped
order. To restore the order $(a,+,x),(b,+,y)$, replace $\zeta$
by a partition $\nu$ sampled according to
\begin{equation}\label{eq:same-sign-up-kernel}
\PP(\nu\mid\alpha,\gamma,\zeta)
=
\frac{U_a(\nu/\alpha;x)U_b(\gamma/\nu;y)}
{\displaystyle
 \sum_{\alpha\subseteq\theta\subseteq\gamma}
 U_a(\theta/\alpha;x)U_b(\gamma/\theta;y)},
\qquad \alpha\subseteq\nu\subseteq\gamma,
\end{equation}
with probability zero outside this range.

For a decreasing pair, consider
\[
(a,-,x),(b,-,y)
\longrightarrow
(b,-,y),(a,-,x),
\]
with fixed endpoints $\gamma\subseteq\alpha$. In this case,
\begin{align}
&\sum_{\gamma\subseteq\nu\subseteq\alpha}
V_a(\alpha/\nu;x)V_b(\nu/\gamma;y)\nonumber\\
&\qquad=
\sum_{\gamma\subseteq\zeta\subseteq\alpha}
V_b(\alpha/\zeta;y)V_a(\zeta/\gamma;x).
\label{eq:same-sign-down}
\end{align}
To restore the order $(a,-,x),(b,-,y)$, replace the current
intermediate partition $\zeta$ by a sample from
\begin{equation}\label{eq:same-sign-down-kernel}
\PP(\nu\mid\alpha,\gamma,\zeta)
=
\frac{V_a(\alpha/\nu;x)V_b(\nu/\gamma;y)}
{\displaystyle
 \sum_{\gamma\subseteq\theta\subseteq\alpha}
 V_a(\alpha/\theta;x)V_b(\theta/\gamma;y)},
\qquad \gamma\subseteq\nu\subseteq\alpha,
\end{equation}
with probability zero outside this range.

Both sampling rules are implemented by enumerating the partitions
between the fixed endpoints, computing the product of the two
transition weights for each candidate, and normalizing these
products to obtain sampling probabilities. Candidates that violate
the strip conditions have zero weight. The normalizing sums are
finite and, by \eqref{eq:same-sign-up} and
\eqref{eq:same-sign-down}, strictly positive whenever the product
of the two current transition weights is positive.

Given the endpoints and the step parameters, these distributions
do not depend on $\zeta$. The size of $\nu$ is not prescribed,
so the powers of $x$ and $y$ in the full edge weights must be
retained. Each commutation identity equates the total local
weights in the two orders. Resampling the intermediate partition
as above therefore gives the conditional distribution for the
restored order while preserving the weights of the fixed exterior.

\subsection{Embedding local identities into a process}

\begin{lemma}[Local transport with a fixed exterior]\label{lem:local-global}
Let $e$ denote the unchanged exterior data, and let $A_e(\nu)$
and $B_e(\kappa)$ be the nonnegative local weights of the original
and reduced models, respectively. Suppose that, for each $e$,
a Markov kernel $T_e$ satisfies
\[
\sum_\kappa B_e(\kappa)T_e(\kappa,\nu)
=
\frac{A_e(\nu)}{f},
\]
where $f>0$ is independent of $e$.
Let $H(e)$ be a common nonnegative exterior weight such that
\[
Z_A=\sum_{e,\nu}H(e)A_e(\nu),
\qquad
Z_B=\sum_{e,\kappa}H(e)B_e(\kappa)
\]
are finite and positive. Then the update that keeps $e$ fixed
and applies $T_e$ to the local variables transports the normalized
measure with weights $H(e)B_e(\kappa)$ to the normalized measure
with weights $H(e)A_e(\nu)$.
\end{lemma}

\begin{proof}
Multiplying the local identity by $H(e)$, summing over $e,\nu$,
and using the stochasticity of $T_e$, we obtain
\[
Z_A
=
f\sum_{e,\kappa}H(e)B_e(\kappa)
  \sum_\nu T_e(\kappa,\nu)
=
fZ_B.
\]
Consequently, for each $(e,\nu)$,
\[
\sum_\kappa
\frac{H(e)B_e(\kappa)}{Z_B}\,T_e(\kappa,\nu)
=
\frac{H(e)A_e(\nu)}{fZ_B}
=
\frac{H(e)A_e(\nu)}{Z_A},
\]
which is the desired transport identity.
Nonnegativity justifies all rearrangements of sums.
\end{proof}

For bulk, boundary, and same-sign moves, the factor $f$ in
Lemma~\ref{lem:local-global} is $C_{ab}(xy)$, $F_{c,a}(hx)$,
and $1$, respectively. These factors depend only on the model parameters, not on the
partitions. Hence, if the current partition sequence has the
process law associated with the labelled word after a local
reduction, the corresponding reverse kernel produces the process
law associated with the labelled word before that reduction.
Applying this argument successively to the recorded moves proves
correctness of the full reverse reconstruction.

\section{A finite commutation--reflection history}\label{sec:history}

A reduction acts only on the labelled word and its parameters; no random partitions are needed to build it. Record each move together with its location, types, and \emph{pre-move} parameters. This information is essential when the history is replayed backwards.

First, repeatedly replace adjacent $+,-$ pairs by $-,+$ pairs until the word is stably sorted, with every minus before every plus. Denote the resulting word by $\bar w$, with its parameters reordered along with the labels. Its sign sequence is
\[
-\,\cdots\,-\,+\,\cdots\,+.
\]
This initial sorting uses at most $N(N-1)/2$ adjacent swaps.

Starting from the sorted labelled word, process each labelled
letter once. A minus letter is moved to the left boundary,
reflected, moved to the right boundary, reflected again, and
returned to its starting slot. A plus letter follows the opposite
boundary order: it is moved first to the right boundary, then to
the left boundary, and finally back to its starting slot. Moving
a labelled letter between slots is done by successive adjacent
swaps.

\begin{lemma}[One reflection period]\label{lem:period}
After a plus or minus letter has visited both boundaries in the
order described above and returned to its starting slot, the
ordering of the labelled letters is the same as before, and the
travelling letter has recovered its original sign. Every
opposite-sign swap encountered during this procedure is a
reduction of the form
\[
+,-\ \longrightarrow\ -,+.
\]
The parameter carried by the travelling letter is multiplied by
$u^2v^2$. Consequently, after all $N$ labelled letters have been
processed once, the sorted labelled word is restored, with every
parameter multiplied by
\begin{equation}\label{eq:r-def}
r=(uv)^2.
\end{equation}
One period uses $2N$ reflections and exactly $2N(N-1)$ adjacent
swaps.
\end{lemma}
\begin{proof}
Consider a minus letter in slot $j$. Since the word is sorted,
moving it to the left boundary crosses only minus letters. The
left reflection changes its sign to $+$ and multiplies its
parameter by $u^2$.

The letter then moves from the left boundary to the right boundary.
Whenever it crosses a letter of the opposite sign, the adjacent
pair is
\[
+,-\ \longrightarrow\ -,+.
\]
All remaining swaps are between two plus letters. At the right
boundary, the second reflection changes the travelling letter
back to $-$ and multiplies its parameter by $v^2$.

On the return from the right boundary to slot $j$, every
opposite-sign swap again has the form
\[
+,-\ \longrightarrow\ -,+,
\]
and all remaining swaps are between two minus letters.

The three parts of the route use
\[
(j-1)+(N-1)+(N-j)=2(N-1)
\]
adjacent swaps. The other labelled letters never swap with one
another, so their relative order is unchanged. When the travelling
letter returns to slot $j$, the full ordering is therefore restored.
The two reflections restore its original sign and multiply its
parameter by $u^2v^2$.

The argument for a plus letter is the same with left and right
interchanged. Hence processing all $N$ labelled letters uses
$2N$ reflections and exactly $2N(N-1)$ adjacent swaps, and
multiplies every parameter by $u^2v^2$.
\end{proof}

The resulting sequence of swaps and reflections is recorded in
Algorithm~\ref{alg:history}.

Before stating the algorithm, we specify how moves are recorded.
A \emph{recorded adjacent swap} consists of the following two
operations:

\begin{enumerate}
\item append to the ordered list $\mathcal H$ the current slot of
the pair and the complete pre-swap data of the two labelled letters,
including their labels, types, signs, and parameters;
\item interchange the two labelled letters.
\end{enumerate}

Moving a labelled letter from one slot to another by recorded
adjacent swaps means repeatedly applying this operation to the
next adjacent pair in the required direction.

A \emph{recorded left reflection} appends to $\mathcal H$ the
pre-reflection data of the labelled letter and the information
that the reflection occurs at the left boundary; it then changes
the sign from $-$ to $+$ and replaces $x_i$ by $u^2x_i$.
A \emph{recorded right reflection} is defined analogously: it
records the pre-reflection data, changes the sign from $+$ to $-$,
and replaces $x_i$ by $v^2x_i$.

\begin{algorithm}[tb]
\caption{Build a reduction history}\label{alg:history}
\begin{algorithmic}[1]
\Require A labelled word $w$ and a period cutoff $K\in\N$.
\State Initialize $\mathcal H$ as an empty ordered list.

\While{an adjacent $+,-$ pair exists}
  \State Perform a recorded adjacent swap on this pair.
\EndWhile
\State Let $\bar w$ be the resulting sorted labelled word.

\For{$k=1,\ldots,K$}
  \For{$j=1,\ldots,N$}
    \State Let $i$ be the label of the letter currently in slot $j$.

    \If{letter $i$ is a minus letter}
      \State Move $i$ to slot $1$ by recorded adjacent swaps.
      \State Apply a recorded left reflection to $i$.
      \State Move $i$ to slot $N$ by recorded adjacent swaps.
      \State Apply a recorded right reflection to $i$.
    \Else
      \State Move $i$ to slot $N$ by recorded adjacent swaps.
      \State Apply a recorded right reflection to $i$.
      \State Move $i$ to slot $1$ by recorded adjacent swaps.
      \State Apply a recorded left reflection to $i$.
    \EndIf

    \State Move $i$ back to slot $j$ by recorded adjacent swaps.
  \EndFor
\EndFor

\State \Return $\mathcal H$ and the final labelled word,
which is $\bar w$ with every parameter multiplied by $r^K$.
\end{algorithmic}
\end{algorithm}

\subsection{The terminal seed and its truncation}\label{subsec:seed}

Put $z=uv$. If every edge parameter is zero, every consecutive pair of partitions must be equal. The common partition has weight
\begin{equation}\label{eq:h-seed}
h_{c_l,c_r}(\lambda)z^{\size{\lambda}},\qquad
h_{c_l,c_r}(\lambda)
=\frac{\beta_{c_l}(\lambda)\beta_{c_r}(\lambda)}{b_\lambda}.
\end{equation}
Thus define
\begin{equation}\label{eq:seed-law}
\seedZ(z)=\sum_{\lambda\in\Y}h_{c_l,c_r}(\lambda)z^{\size{\lambda}},
\qquad
\rho(\lambda)=\frac{h_{c_l,c_r}(\lambda)z^{\size{\lambda}}}{\seedZ(z)}.
\end{equation}
The scalar $\seedZ$ is the terminal normalization in
\cite[(3.2)]{LiShi2026}. In particular, the summation is over all allowed partitions, not just single rows.

For $M\in\N$, let
\begin{equation}\label{eq:seed-trunc}
\seedZ_M(z)=\sum_{\size{\lambda}\leq M}h_{c_l,c_r}(\lambda)z^{\size{\lambda}},
\qquad
\rho_M(\lambda)=\frac{\ind_{\{\size{\lambda}\leq M\}}
 h_{c_l,c_r}(\lambda)z^{\size{\lambda}}}{\seedZ_M(z)}.
\end{equation}
This is a finite categorical distribution with weights given by finite box products. It can be sampled without evaluating $\seedZ(z)$. When used to initialize a path, $\rho_M$ means the law of
$(\Lambda,\ldots,\Lambda)$ with $\Lambda\sim\rho_M$.

Proposition~\ref{prop:exact-seed} gives an exact sampler for $\rho$
without evaluating $\mathfrak Z(z)$; Algorithm~2 uses the truncated
seed $\rho_M$.

\subsection{Reverse replay}

\begin{algorithm}[tb]
\caption{Approximate doubly free-boundary Macdonald sampler}
\label{alg:sampler}
\begin{algorithmic}[1]
\Require Word $w$, boundary data $(u,v,c_l,c_r)$, and cutoffs $K,M\in\N$.
\State Construct the history $\mathcal H$ using Algorithm~\ref{alg:history}.
\State Sample $\Lambda\sim\rho_M$ and set
$\lambda^{(i)}\gets\Lambda$ for $i=0,\ldots,N$.
\For{recorded moves in $\mathcal H$, in reverse order}
  \State Read the move location and pre-move letter data from the record.
  \If{the recorded move was an adjacent swap $(+,-)\to(-,+)$}
    \State Read the outer partitions $\lambda,\mu$ and current valley $\kappa$.
    \State Draw $G\sim G_{ab}(xy)$.
    \State Set $m\gets\size{\lambda}+\size{\mu}-\size{\kappa}+G$.
    \State Replace the valley partition $\kappa$
    by a sample from \eqref{eq:bulk-shape}.
  \ElsIf{the recorded move was a boundary reflection}
    \State Read the current boundary partition $\kappa$ and its neighbor $\eta$.
    \State Set $(h,c)\gets(u,c_l)$ on the left, or $(v,c_r)$ on the right.
    \State Draw $B\sim B_{c,a}(hx)$ and set
    $m\gets2\size{\eta}-\size{\kappa}+B$.
    \State Replace the boundary partition $\kappa$
    by a sample from \eqref{eq:boundary-shape}.
  \Else \Comment{Same-sign adjacent swap}
    \State Resample the middle partition using
    \Statex \hspace{\algorithmicindent}\hspace{\algorithmicindent}
    \eqref{eq:same-sign-up-kernel} for $(+,+)$ or
    \eqref{eq:same-sign-down-kernel} for $(-,-)$.
  \EndIf
\EndFor
\State \Return the partition sequence on the original word $w$.
\end{algorithmic}
\end{algorithm}

The seed and all scalar increments are drawn using independent randomness.
The scalar distributions depend only on the recorded types and parameters,
so these increments may be sampled before reverse replay. Each shape draw
uses the current partition path and the recorded pre-move weights; it
cannot be determined from the history alone.

Let $n_+$ and $n_-$ be the numbers of plus and minus letters, and define
\[
  I(w)=\#\{(i,j):1\le i<j\le N,\ \sigma_i=+,\ \sigma_j=-\}.
\]
Initial sorting uses exactly $I(w)$ swaps. Since $n_++n_-=N$,
$I(w)\le n_+n_-\le\lfloor N^2/4\rfloor$.
In one reflection period, each travelling letter crosses every other
letter twice, once with the same sign and once with the opposite sign,
and undergoes two boundary reflections. Thus the history contains
\[
  L_{\mathrm{bulk}}=I(w)+KN(N-1),\qquad
  L_{\mathrm{bdry}}=2KN,\qquad
  L_{\mathrm{same}}=KN(N-1)
\]
bulk, boundary, and same-sign moves, respectively. Consequently,
\begin{equation}\label{eq:move-count}
  L_K=I(w)+2KN^2
  \le n_+n_-+2KN^2
  \le\left\lfloor\frac{N^2}{4}\right\rfloor+2KN^2.
\end{equation}
This is the exact number of reverse local updates; $M$ does not affect it.
At unit cost per letter or scalar-parameter operation, the history can
be constructed in $O(N+L_K)$ operations. The word and its records occupy
$O(N+L_K)$ entries. Initial sorting can be implemented by scanning the
word and moving each minus letter left past the preceding plus letters.

Let $n_L$ and $n_R$ be the numbers of letters of types $L$ and $R$,
so that $n_L+n_R=N$. Define
\[
  I_{\mathrm{same}}
  :=\#\{(i,j):1\le i<j\le N,\ \sigma_i=+,\ \sigma_j=-,\ a_i=a_j\},
  \qquad
  I_{\mathrm{mix}}:=I(w)-I_{\mathrm{same}}.
\]
Initial sorting requires $I_{\mathrm{same}}$ calls to the rejection
sampler of Lemma~3.1 and $I_{\mathrm{mix}}$ direct Bernoulli draws.
Each reflection period has $n_L(n_L-1)+n_R(n_R-1)$ same-type bulk
moves, each requiring one call to Lemma~3.1, and $2n_Ln_R$ mixed-type
bulk moves, each requiring only a direct Bernoulli draw.

Every labelled letter visits each boundary once per period.
For a boundary of type $c$, let $d(c)$ count the zero increments
and $b(c)$ count the additional Bernoulli draws in one period.
The boundary increment laws give
\[
\begin{array}{c|cccc}
  c & \mathrm{el} & \mathrm{oa} & \mathrm{deel} & \mathrm{eoa} \\
  \hline
  d(c) & 0 & 0 & n_L & n_R \\
  b(c) & n_R & n_L & 0 & 0
\end{array}
\]
The remaining $N-d(c)$ boundary increments each require one call
to Lemma~3.1.

The total number of bulk and boundary moves is
\[
  D=L_{\mathrm{bulk}}+L_{\mathrm{bdry}}
   =I(w)+KN(N+1).
\]
The exact number of calls to Lemma~3.1 is
\[
\begin{aligned}
  n_{\mathrm{NB}}
  &=I_{\mathrm{same}}
    +K\bigl[n_L^2+n_R^2+N-d(c_l)-d(c_r)\bigr] \\
  &=D-I_{\mathrm{mix}}-2Kn_Ln_R-K\bigl[d(c_l)+d(c_r)\bigr].
\end{aligned}
\]
Since
\[
  I_{\mathrm{same}}
  \le \left\lfloor\frac{n_L^2}{4}\right\rfloor
     +\left\lfloor\frac{n_R^2}{4}\right\rfloor,
\]
we obtain the sharper bound
\[
  n_{\mathrm{NB}}
  \le \left\lfloor\frac{n_L^2}{4}\right\rfloor
     +\left\lfloor\frac{n_R^2}{4}\right\rfloor
     +K\bigl[n_L^2+n_R^2+N-d(c_l)-d(c_r)\bigr].
\]
For fixed type counts and boundary types, this bound is attained by
balancing the signs within each type and placing all plus letters
before all minus letters.

The exact number of direct Bernoulli draws is
\[
  n_{\mathrm{Bern}}
  =I_{\mathrm{mix}}+2Kn_Ln_R+K\bigl[b(c_l)+b(c_r)\bigr]
  \le (2K+1)n_Ln_R+K\bigl[b(c_l)+b(c_r)\bigr].
\]
This count excludes the acceptance tests inside Lemma~3.1.
The $K[d(c_l)+d(c_r)]$ deterministic increments are set to zero
without random sampling.

Index the calls to Lemma~3.1 by $j=1,\ldots,n_{\mathrm{NB}}$.
Call $j$ samples $\mathrm{NB}_{Q_j}(A_j,s_j)$, where
\[
  (A_j,Q_j)\in\{(t,q),(q,t),(qt,q^2),(qt,t^2)\}.
\]
The argument $s_j$ is the corresponding recorded value of $xy$, $hx$,
or $(hx)^2$, where $h=u$ for a left boundary reflection and $h=v$
for a right boundary reflection. Recall that
\[
\begin{aligned}
  H_Q&=(1-Q)^{-\lceil(1-Q)^{-1}\rceil},\\
  K_{A,Q}(s)&=\prod_{r=0}^{\infty}
  \frac{1-AsQ^r}{1-sQ^r}.
\end{aligned}
\]
The coefficient bound in Lemma~3.1 and the first factor of this product give
\[
  \frac{1-As}{1-s}\le K_{A,Q}(s)\le\frac{H_Q}{1-s}.
\]
Each trial in call $j$ succeeds independently with probability
\[
  p_j=\frac{(1-s_j)K_{A_j,Q_j}(s_j)}{H_{Q_j}}.
\]
Let $T_j$ count its trials, including the successful trial. Then $T_j$
has a geometric distribution on $\{1,2,\ldots\}$, with
\[
\begin{aligned}
  \mu_j&:=\mathbb E[T_j]
  =\frac{H_{Q_j}}{(1-s_j)K_{A_j,Q_j}(s_j)},\\
  \operatorname{Var}(T_j)&=\mu_j(\mu_j-1).
\end{aligned}
\]
Writing $U_j=H_{Q_j}/(1-A_js_j)$, the preceding lower bound implies
$\mu_j\le U_j\le H_{Q_j}/(1-A_j)$.
For $\vartheta=\max(q,t)$, all calls satisfy $A_j,Q_j\le\vartheta$.
Since $H_Q$ is nondecreasing in $Q$, set
\[
  U_*:=\frac{H_\vartheta}{1-\vartheta},\qquad U_j\le U_*.
\]
The calls use independent randomness. Hence their total trial count
$T_{\mathrm{tot}}=\sum_jT_j$ satisfies
\[
\begin{aligned}
  \mathbb E[T_{\mathrm{tot}}]
    &=\sum_j\mu_j
      \le\sum_jU_j\le n_{\mathrm{NB}}U_*\le DU_*,\\
  \operatorname{Var}(T_{\mathrm{tot}})
    &=\sum_j\mu_j(\mu_j-1)
      \le\sum_jU_j(U_j-1)
      \le DU_*(U_*-1).
\end{aligned}
\]

A trial proposing $J$ computes $k_J(A_j,Q_j)$ using $J$
finite-product factors. Including the acceptance test, this takes
$O(J+1)$ operations under unit-cost arithmetic and random sampling
primitives. Define
\[
  x_*:=\max_i x_i,\qquad b_*:=\max(u,v)x_*,\qquad
  s_*:=\max(x_*^2,b_*)<1.
\]
Recorded letter parameters never exceed $x_*$, so every call satisfies
$s_j\le s_*$. In the $j$th run of the rejection sampler in Lemma~3.1,
let $J_{j,r}$ be the candidate integer drawn in the $r$th attempt.
Its geometric distribution gives
$\mathbb E[J_{j,r}]=s_j/(1-s_j)$.
The event that trial $r$ is reached depends only on previous trials,
so it is independent of $J_{j,r}$. Therefore
\[
\begin{aligned}
  \mathbb E\!\left[
    \sum_{j=1}^{n_{\mathrm{NB}}}\sum_{r=1}^{T_j}(1+J_{j,r})
  \right]
  &=\sum_{j=1}^{n_{\mathrm{NB}}}\frac{\mu_j}{1-s_j} \\
  &\le \frac{n_{\mathrm{NB}}U_*}{1-s_*}
   \le \frac{DU_*}{1-s_*}.
\end{aligned}
\]
Direct Bernoulli draws and deterministic assignments require $O(D)$
additional operations. Thus the expected arithmetic cost of sampling
all scalar increments is
\[
  O\!\left(D+\frac{n_{\mathrm{NB}}U_*}{1-s_*}\right),
\]
and is at most $O(DU_* /(1-s_*))$.

Let $p(n)$ be the number of partitions of $n$.
If $\Lambda\sim\rho_M$ is sampled by enumerating all partitions
of size at most $M$, the number of candidates is
\[
  C_{\mathrm{seed}}=\sum_{n=0}^{M}p(n).
\]

The history has finitely many records, every scalar sampler terminates
almost surely, and every shape candidate set is finite. Hence
Algorithm~\ref{alg:sampler} terminates almost surely.
Bound~\eqref{eq:move-count} controls the number of local updates,
but shape enumeration can be expensive, so it does not imply a
polynomial running time.

\section{Total variation convergence}\label{sec:TV}

Let $\bar w$ be the labelled word obtained by the initial sorting
in Section~\ref{sec:history}, with all minus letters preceding
all plus letters. For $s\in(0,1]$, let $\xi_s$ be its normalized process after multiplying \emph{every} edge parameter by $s$, and let $\ZZ_s$ be its partition function. At $s=0$, let $\xi_0$ be the law of
$(\Lambda,\ldots,\Lambda)\in\mathbb Y^{N+1}$, where $\Lambda\sim\rho$. Thus
\begin{equation}\label{eq:Z0}
\ZZ_0=\seedZ(z),\qquad z=uv.
\end{equation}
Starting from a partition sequence on the final word produced by
Algorithm~\ref{alg:history}, apply the local updates in reverse order:
first undo the $K$ reflection periods, then the initial sorting.
Let $\mathcal R_K(\lambda,\mu)$ be the probability that this procedure
outputs the partition sequence $\mu$ when started from $\lambda$. For an initial law $\alpha$, write $\alpha\mathcal R_K$ for the
output law after these updates:
\[
  (\alpha\mathcal R_K)(\mu)
  =\sum_{\lambda}\alpha(\lambda)\mathcal R_K(\lambda,\mu),
\]
where the sum is over admissible partition sequences on the final word. Lemma~\ref{lem:local-global}, applied at each move, gives the exact identity
\begin{equation}\label{eq:exact-replay}
\xi_{r^K}\mathcal R_K=\pi_w,
\qquad r=z^2.
\end{equation}
The actual output law is
\begin{equation}\label{eq:output-law}
\widehat\pi_{K,M}=\rho_M\mathcal R_K,
\end{equation}
where $\rho_M$ is embedded on constant paths as above.

We use the convention
\[
\tv{\mu-\nu}=\sup_A|\mu(A)-\nu(A)|
=\frac12\sum_x|\mu(x)-\nu(x)|
\]
on a countable state space.

\begin{theorem}[Two-cutoff approximation]\label{thm:TV}
For the process \eqref{eq:target-law} under \eqref{eq:regime}, Algorithm~\ref{alg:sampler} is well defined and satisfies
\begin{equation}\label{eq:TV-exact}
\boxed{\quad
\tv{\widehat\pi_{K,M}-\pi_w}
\leq1-\frac{\seedZ_M(z)}{\ZZ_{r^K}}.
\quad}
\end{equation}
In particular, with $C=\ZZ_1/\ZZ_0-1<\infty$,
\begin{equation}\label{eq:TV-rate}
\tv{\widehat\pi_{K,M}-\pi_w}
\leq Cr^K+\rho(\size{\Lambda}>M).
\end{equation}
Consequently,
\begin{equation}\label{eq:TV-limit}
\lim_{\substack{K\to\infty\\ M\to\infty}}
\tv{\widehat\pi_{K,M}-\pi_w}=0.
\end{equation}
The limit holds along every sequence for which both cutoffs diverge.
\end{theorem}

\begin{proof}
For a partition sequence
$\boldsymbol\lambda=(\lambda^{(0)},\ldots,\lambda^{(N)})
\in\mathbb{Y}^{N+1}$,  define its total skew degree
\begin{equation}\label{eq:degree}
D(\boldsymbol\lambda)=\sum_{i=1}^N
\left|\size{\lambda^{(i)}}-\size{\lambda^{(i-1)}}\right|.
\end{equation}
By homogeneity, its weight under uniform edge scaling is
\begin{equation}\label{eq:scale-weight}
W_s(\boldsymbol\lambda)=s^{D(\boldsymbol\lambda)}W_1(\boldsymbol\lambda).
\end{equation}
On the allowed support, consecutive partitions differ by adding or
removing a horizontal or vertical strip.
If $D(\boldsymbol\lambda)=0$, every consecutive pair has the same size,
so each such strip is empty and
$\lambda^{(0)}=\cdots=\lambda^{(N)}$.
The converse is immediate. Therefore the $D=0$ contribution is exactly $\seedZ(z)$. Since $\ZZ_1<\infty$, dominated convergence gives
\begin{equation}\label{eq:Z-continuity}
\ZZ_s\longrightarrow\ZZ_0\qquad(s\downarrow0).
\end{equation}
Moreover, for $0\leq s\leq1$,
\begin{equation}\label{eq:Z-linear}
0\leq\ZZ_s-\ZZ_0
=\sum_{D(\boldsymbol\lambda)\geq1}s^{D(\boldsymbol\lambda)}W_1(\boldsymbol\lambda)
\leq s(\ZZ_1-\ZZ_0).
\end{equation}

Let $E_M$ be the event that every partition is equal and their common size is at most $M$. Under $\xi_s$,
\begin{equation}\label{eq:event-prob}
\xi_s(E_M)=\frac{\seedZ_M(z)}{\ZZ_s},\qquad
\xi_s(\,\cdot\mid E_M)=\rho_M.
\end{equation}
For any probability measure $\mu$ and event $E$ of positive probability,
$\tv{\mu-\mu(\cdot\mid E)}=1-\mu(E)$. This follows by splitting the half-sum defining total variation into $E$ and its complement. Hence
\begin{equation}\label{eq:initial-TV}
\tv{\xi_s-\rho_M}=1-\frac{\seedZ_M(z)}{\ZZ_s}.
\end{equation}
Both $\rho_M$ and $\xi_{r^K}$ are processed using the same reverse
updates $\mathcal R_K$. The total variation error after these updates
is at most the error before them; see
\cite[Proposition~3(d)]{RobertsRosenthal2004}.
Using \eqref{eq:exact-replay}, \eqref{eq:output-law}, and
\eqref{eq:initial-TV}, we obtain
\[
\begin{aligned}
  \|\widehat\pi_{K,M}-\pi_w\|_{\mathrm{TV}}
  &=\|\rho_M\mathcal R_K-\xi_{r^K}\mathcal R_K\|_{\mathrm{TV}}\\
  &\le\|\rho_M-\xi_{r^K}\|_{\mathrm{TV}}\\
  &=1-\frac{\mathfrak Z_M(z)}{\ZZ_{r^K}}.
\end{aligned}
\]
This proves \eqref{eq:TV-exact}.

Finally,
\begin{align*}
1-\frac{\seedZ_M(z)}{\ZZ_s}
&=\frac{\ZZ_s-\ZZ_0}{\ZZ_s}
 +\frac{\ZZ_0-\seedZ_M(z)}{\ZZ_s}\\
&\leq s\left(\frac{\ZZ_1}{\ZZ_0}-1\right)
 +\left(1-\frac{\seedZ_M(z)}{\seedZ(z)}\right).
\end{align*}
The final parenthesis is $\rho(\size{\Lambda}>M)$. This proves \eqref{eq:TV-rate}; the seed tail tends to zero because $\rho$ is a probability measure. Almost-sure termination was established after Algorithm~\ref{alg:sampler}.
\end{proof}

\begin{corollary}[Exact seed and observables]\label{cor:observables}
If the seed $\rho$ can be sampled exactly, then the seed-cutoff error is absent and
\[
\tv{\rho\mathcal R_K-\pi_w}
\leq1-\frac{\seedZ(z)}{\ZZ_{r^K}}\leq Cr^K.
\]
\end{corollary}

\section{Explicit error bounds and implementation}
\label{sec:explicit-bounds}
The bound in Theorem~\ref{thm:TV} involves a ratio of partition functions.
We give explicit truncation bounds and an expected cost bound for
direct enumeration.

\subsection{Reflection depth and enumeration cost}
\label{sec:reflection-cost}
Set
\begin{equation}\label{eq:explicit-bound-parameters}
  \vartheta=\max(q,t),\qquad
  L_*=(1-\vartheta)^{-1},\qquad
  x_* = \max_i x_i,\qquad
  b_* = \max(u,v)x_*.
\end{equation}
Define
\begin{equation}\label{eq:explicit-bound-constants}
\begin{aligned}
  A_*&=2N(1+L_*)\frac{b_*}{1-b_*},&
  B_*&=2N(N-1)L_*\frac{x_*^2}{1-x_*^2},\\
  H(s)&=\frac{A_*s}{1-r}+\frac{B_*s^2}{1-r^2},&
  r&=(uv)^2.
\end{aligned}
\end{equation}

\begin{proposition}[Reflection-depth bound]
\label{prop:reflection-depth}
For $0\le s\le1$,
\begin{equation}\label{eq:reflection-depth-bound}
  \log\frac{\ZZ_s}{\ZZ_0}\le H(s),\qquad
  1-\frac{\ZZ_0}{\ZZ_s}\le1-e^{-H(s)}\le H(s).
\end{equation}
\end{proposition}
\begin{proof}
For every factor $K_{A,Q}$ used in the history, $Q\le\vartheta$.
Thus, for $0\le y<1$,
\[
  \log K_{A,Q}(y)
  =\sum_{n\ge1}\frac{1-A^n}{1-Q^n}\frac{y^n}{n}
  \le L_*\frac{y}{1-y}.
\]
Together with $\log(1+y)\le y$, the bulk and boundary factor formulas
\eqref{eq:C-factors}, \eqref{eq:F-table}, and \eqref{eq:F-dual} give
\[
  \log C_{ab}(y)\le L_*\frac{y}{1-y},\qquad
  \log F_{c,a}(y)\le(1+L_*)\frac{y}{1-y}.
\]
Record the pre-move data in one reflection period started from
$\bar w$ with edge parameters $x_1,\ldots,x_N$. Define
\[
  F(s):=
  \prod_{\text{opposite-sign swaps}} C_{ab}(s^2xy)
  \prod_{\text{boundary reflections}} F_{c,a}(shx),
\]
where $x,y$ are the recorded pre-move edge parameters, and
$h=u$ at the left boundary or $h=v$ at the right boundary.
Multiplying every initial edge parameter by $s$ multiplies
every recorded edge parameter by $s$.
Each factor uses the types and parameters immediately before the
corresponding move, with $h=u$ at the left boundary and $h=v$ at
the right boundary. Same-sign swaps contribute a factor of one.
Throughout the period, bulk arguments are at most $s^2x_*^2$ and
boundary arguments at most $sb_*$. The move counts give
$\log F(s)\le A_*s+B_*s^2$.
The reduction identities give $\ZZ_s=F(s)\ZZ_{rs}$.
Iterating and using $\ZZ_{sr^j}\to\ZZ_0$ yields
\[
  \log\frac{\ZZ_s}{\ZZ_0}
  =\sum_{j\ge0}\log F(sr^j)
  \le\frac{A_*s}{1-r}+\frac{B_*s^2}{1-r^2}.
\]
The remaining inequalities follow by exponentiation.
\end{proof}

\begin{theorem}[Expected cost of direct enumeration]
\label{thm:expected-running-time}
Fix $q,t,u,v,x_0\in(0,1)$ and suppose that $0<x_i\le x_0$ for every letter.
Let $N\ge1$ and $K,M\in\mathbb Z_{\ge0}$.
Implement Algorithm~\ref{alg:sampler} by direct enumeration for the
initial partition and all conditional shape draws, using the weights in
Appendix~A. Count arithmetic operations and geometric and Bernoulli draws
at unit cost. If $\mathcal T_{K,M}$ is the total operation count, including
history construction and initial sampling, then
\[
  \mathbb E[\mathcal T_{K,M}]
  \le C_0\left[(K+1)(N+1)^2e^{c_0N}
    +(M+1)^3e^{c_p\sqrt M}\right],
  \qquad c_p=\pi\sqrt{2/3},
\]
where $C_0,c_0>0$ depend only on $q,t,u,v,x_0$, uniformly over the word
and boundary types.
\end{theorem}
\begin{proof}
Choose $\delta>0$ with $2\delta<\min\{-\log u,-\log x_0\}$.
For any word encountered during replay, replace $u$ by $ue^\delta$
and multiply every edge parameter by $e^\delta$.
The partition-function estimate in Appendix~B, applied to these
modified parameters, gives
\[
  \widetilde Z\le A e^{BN^2},
\]
with $A,B>0$ depending only on the fixed parameters.

At any stage of the reverse reconstruction, let
$\nu(\boldsymbol\lambda)$ be the probability that the sampler has
partition sequence $\boldsymbol\lambda$, and let
$W(\boldsymbol\lambda)$ be its unnormalized weight for the current word.
Then
\[
  \nu(\boldsymbol\lambda)\le W(\boldsymbol\lambda).
\]
Indeed, initially, a constant sequence $(\Lambda,\ldots,\Lambda)$
with $\size{\Lambda}\le M$ has probability equal to its terminal
weight divided by $\mathfrak Z_M(uv)\ge1$; all other sequences have
probability zero. Each reverse update transports the current
unnormalized weight to the next unnormalized weight divided by a
scalar factor at least one. Thus the inequality is preserved
throughout the reconstruction.

Put $S=\max_{0\le i\le N}\size{\lambda^{(i)}}$. The total size
variation $D(\boldsymbol\lambda)$ defined in \eqref{eq:degree} gives
\[
  S\le \size{\lambda^{(0)}}+D(\boldsymbol\lambda).
\]
Replacing $u$ by $ue^\delta$ and every edge parameter $x$ by
$xe^\delta$ multiplies the weight by
$e^{\delta(\size{\lambda^{(0)}}+D(\boldsymbol\lambda))}$.
Consequently,
\[
  \widetilde W(\boldsymbol\lambda)
  =W(\boldsymbol\lambda)
    e^{\delta(\size{\lambda^{(0)}}+D(\boldsymbol\lambda))}
  \ge W(\boldsymbol\lambda)e^{\delta S}.
\]
Therefore
\[
  \mathbb E_\nu[e^{\delta S}]
  =\sum_{\boldsymbol\lambda}
     \nu(\boldsymbol\lambda)e^{\delta S}
  \le\sum_{\boldsymbol\lambda}\widetilde W(\boldsymbol\lambda)
  =\widetilde Z
  \le Ae^{BN^2}.
\]

Let $S$ be the largest partition size after the update.
For a bulk or boundary update, every candidate has the prescribed
size $m\le S$. For a same-sign update, every candidate has size
at most the larger of the two fixed outer partition sizes,
which is also at most $S$. By \eqref{eq:partition-count},
the number of candidates is therefore at most
\[
  \sum_{n=0}^{S}p(n)\le(S+1)e^{c_p\sqrt S}.
\]
Although $S$ is random, the exponential-moment bound proved above
holds at every reconstruction stage.
Using \eqref{eq:branching-products} and
\eqref{eq:coeff-implementation}, each candidate weight involves
$O(S+1)$ box factors. Each factor contains a fixed number of powers
with exponents at most $S$. Evaluating these powers by repeated
multiplication takes $O(S+1)$ operations per factor, giving
$O((S+1)^2)$ operations per candidate. Thus one shape draw costs at most
$C_1(S+1)^3e^{c_p\sqrt S}$, including candidate generation,
normalization, and sampling.
Put $a=c_p+1$. Since $(S+1)^3\le C_2e^{\sqrt S}$ and
\[
  a\sqrt S\le\frac{\delta S}{N+1}
               +\frac{a^2(N+1)}{4\delta},
\]
H\"older's inequality gives
\[
  \mathbb E_\nu[(S+1)^3e^{c_p\sqrt S}]
  \le C_2e^{a^2(N+1)/(4\delta)}
      \bigl(\mathbb E_\nu[e^{\delta S}]\bigr)^{1/(N+1)}
  \le C_3e^{c_0N}.
\]
Summing over the $L_K$ updates gives an expected shape cost at most
$C_4L_Ke^{c_0N}$.
The initial enumeration costs at most
$C_5(M+1)^3e^{c_p\sqrt M}$.
History construction, initialization by references to the common
partition, and scalar sampling cost $O(N+L_K)$ in expectation, as shown
in Section~5.2. Finally,
$L_K\le\lfloor N^2/4\rfloor+2KN^2$ proves the claim.
\end{proof}

\subsection{Initial-partition tail and exact sampling}
\label{sec:seed-tail-exact}
Set
\begin{equation}\label{eq:seed-envelope-constants}
  \gamma(q,t)=
  \frac{\sqrt2\,|q-t|}
       {(1-\max(q,t))(1-q)(1-t)},\qquad
  \Gamma=\gamma(q,t)+\pi\sqrt{2/3}.
\end{equation}
With $h_{c_l,c_r}$ defined in \eqref{eq:h-seed},
\eqref{eq:h-bound} and \eqref{eq:partition-count} give,
uniformly over the boundary types,
\begin{equation}\label{eq:seed-coefficient-bound}
  a_n:=\sum_{\size{\lambda}=n}h_{c_l,c_r}(\lambda)
  \le p(n)e^{\gamma(q,t)\sqrt n}
  \le e^{\Gamma\sqrt n},\qquad n\ge0.
\end{equation}
For $z=uv<\eta<1$, put
\[
  C_\eta=\exp\!\left(\frac{\Gamma^2}{4\log(\eta/z)}\right).
\]
The inequality $\Gamma\sqrt n-\delta n\le\Gamma^2/(4\delta)$,
with $\delta=\log(\eta/z)$, gives the common bound
$a_n(z/\eta)^n\le C_\eta$.
Since $\seedZ(z)\ge1$ and $a_nz^n\le C_\eta\eta^n$,
\eqref{eq:seed-law} gives
\begin{equation}\label{eq:seed-tail-bound}
\begin{aligned}
  \rho(\size{\Lambda}>M)
  &=\frac{1}{\seedZ(z)}
    \sum_{n=M+1}^{\infty}a_nz^n\\
  &\le C_\eta\sum_{n=M+1}^{\infty}\eta^n\\
  &=T_M(z,\eta):=
    C_\eta\frac{\eta^{M+1}}{1-\eta}.
\end{aligned}
\end{equation}
Since $0<z<1$, the choice $\eta=\sqrt z$ satisfies $z<\eta<1$.

\begin{corollary}[Explicit truncation bound]
\label{cor:explicit-truncation}
For every $z<\eta<1$,
\begin{equation}\label{eq:parameter-error-bound}
  \|\widehat\pi_{K,M}-\pi_w\|_{\mathrm{TV}}
  \le\min\left\{1,\;1-e^{-H(r^K)}+T_M(z,\eta)\right\}.
\end{equation}
In particular, $H(r^K)\le\varepsilon/2$ and
$T_M(z,\eta)\le\varepsilon/2$ ensure error at most $\varepsilon$.
\end{corollary}
\begin{proof}
Use $1-\mathfrak Z_M(z)/\ZZ_s
\le1-\ZZ_0/\ZZ_s+\rho(\size{\Lambda}>M)$ in
Theorem~\ref{thm:TV}, and apply
\eqref{eq:reflection-depth-bound} and \eqref{eq:seed-tail-bound}
at $s=r^K$.
\end{proof}

\begin{proposition}[Exact initial sampling without a size cutoff]
\label{prop:exact-seed}
The law $\rho$ can be sampled without a size cutoff or evaluation of
$\mathfrak Z(z)$. Fix $z<\eta<1$ and repeat independent trials:
\begin{enumerate}
  \item Draw $J$ with $\mathbb P(J=n)=(1-\eta)\eta^n$, $n\ge0$.
  \item Enumerate the partitions of size $J$ and compute $a_J$.
  If $a_J=0$, restart. Otherwise accept $J$ with probability
  \[
    \alpha_J=\frac{a_J}{C_\eta}\left(\frac z\eta\right)^J.
  \]
  If rejected, restart.
  \item Given acceptance, return a partition $\Lambda$ of size $J$
  with probability $h_{c_l,c_r}(\Lambda)/a_J$.
\end{enumerate}
The procedure terminates almost surely and returns a sample from $\rho$.
\end{proposition}
\begin{proof}
The common coefficient bound gives $0\le\alpha_J\le1$.
For a positive-weight partition $\lambda$ of size $n$, one trial returns
$\lambda$ with probability
\[
  (1-\eta)\eta^n
  \frac{a_n}{C_\eta}\left(\frac z\eta\right)^n
  \frac{h_{c_l,c_r}(\lambda)}{a_n}
  =\frac{1-\eta}{C_\eta}h_{c_l,c_r}(\lambda)z^n.
\]
Thus the success probability is
$(1-\eta)\mathfrak Z(z)/C_\eta>0$, and conditioning on success gives
$\rho$. Each trial enumerates a finite set, proving almost-sure termination.
\end{proof}
Replacing the initial draw in Algorithm~\ref{alg:sampler} by this
procedure gives the output law $\rho\mathcal R_K$ and the bound
\[
  \|\rho\mathcal R_K-\pi_w\|_{\mathrm{TV}}
  \le1-e^{-H(r^K)}.
\]

.

\subsection{A numerical example with fixed Macdonald parameters}
\label{subsec:fixed-qt-example}

\begin{figure}[p]
\centering
\includegraphics[width=\linewidth]{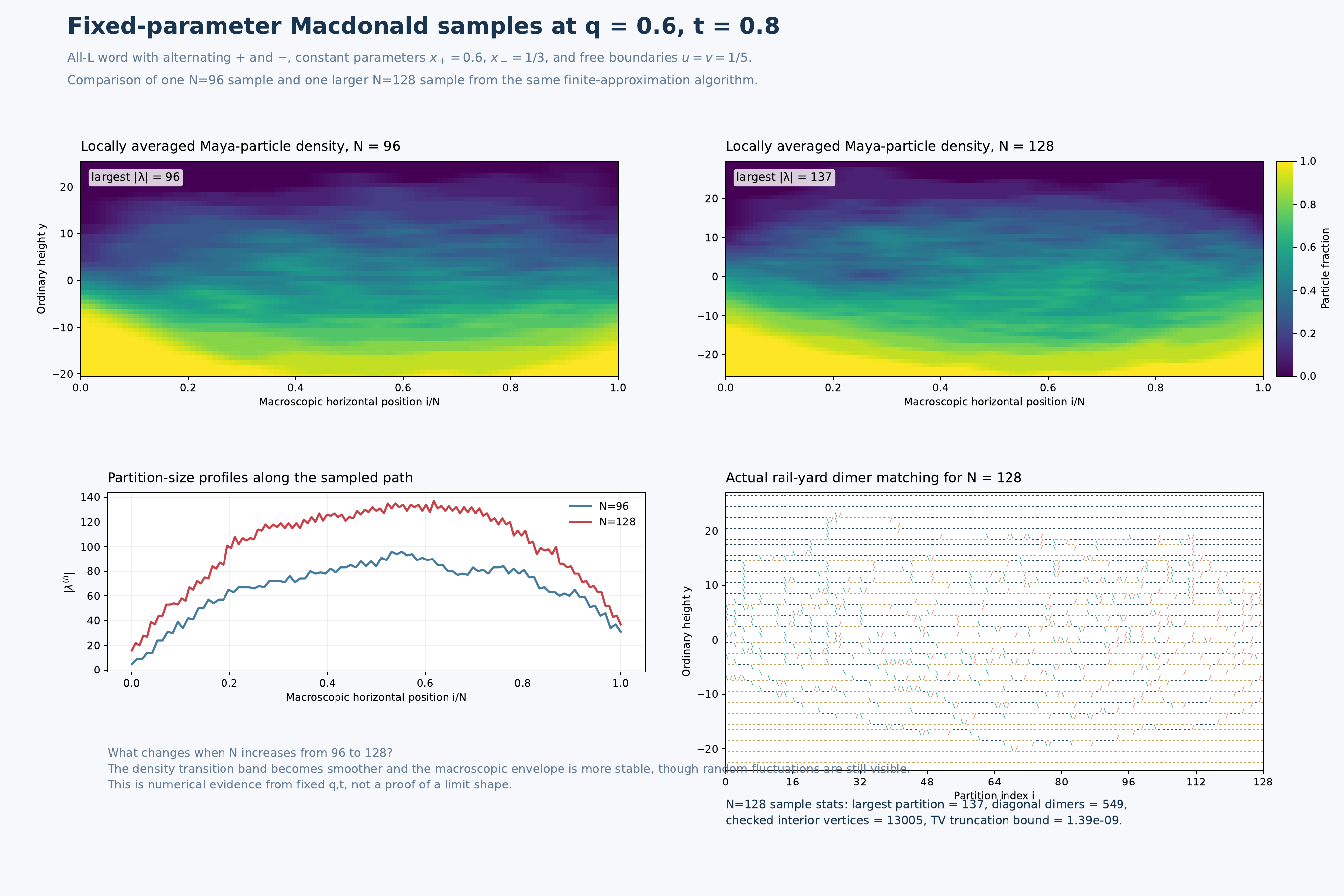}
\caption{An example of Algorithm~\ref{alg:sampler} with fixed
$q=0.6$, $t=0.8$, $u=v=1/5$, and $c_l=c_r=\el$.
The word repeats $(L,+,3/5),(L,-,1/3)$, with $N=96$ or $128$,
and both runs use $K=4$, $M=20$.
Top: locally averaged particle occupancies of the two samples.
Bottom left: the unscaled partition sizes $\size{\lambda^{(i)}}$.
Bottom right: the $N=128$ dimer configuration in ordinary geometric
coordinates.  Each size is represented by one sample;
this visual comparison does not establish a limit shape.}
\label{fig:fixed-qt-example}
\end{figure}

We illustrate Algorithm~\ref{alg:sampler} with the two samples in
Figure~\ref{fig:fixed-qt-example}, on rail-yard graphs with
$N=96$ and $N=128$ columns, respectively.
The Macdonald parameters $q,t$ remain fixed as $N$ increases.

\paragraph{Model and parameter choice.}
We use the rail-yard graph conventions of
\cite[Section~2]{LiShi2026}. Here $N$ is the number of elementary
columns. Column $i$ has type $(a_i,\sigma_i)$, and $x_i$ is the
parameter assigned to its diagonal edges. The partitions
$\lambda^{(i-1)}$ and $\lambda^{(i)}$ encode the boundary profiles
on the two sides of this column. Thus $N$ columns correspond to
a sequence of $N+1$ partitions. We sample this sequence with the
Macdonald weights in \eqref{eq:target-weight}.

For each of these two even values of $N$, set
\begin{equation}\label{eq:fixed-qt-parameters}
q=\frac35,\qquad t=\frac45,\qquad
u=v=\frac15,\qquad c_l=c_r=\el.
\end{equation}
All transitions are of type $L$, with alternating signs and parameters
\begin{equation}\label{eq:fixed-qt-word}
\begin{aligned}
(a_{2k+1},\sigma_{2k+1},x_{2k+1})&=(L,+,3/5),\\
(a_{2k+2},\sigma_{2k+2},x_{2k+2})&=(L,-,1/3),
\qquad 0\leq k<N/2.
\end{aligned}
\end{equation}
Thus the word $w_N$ repeats the pair $(L,+,3/5),(L,-,1/3)$.
Each pair first adds a horizontal strip and then removes a horizontal
strip.  This choice keeps the transition pattern simple while retaining
the non-Schur weights $q\ne t$.  The two $\el$ boundaries impose no parity
restriction: the endpoint partitions are random, not fixed to
$\varnothing$, and retain the weights $\tau_\el$ and $\beta_\el$ in
\eqref{eq:target-weight}.  All parameters satisfy \eqref{eq:regime}; in particular,
$x_{2k+1}x_{2k+2}=1/5<1$.  These are illustrative choices, not parameters
required by the algorithm.  Their rational values permit finite box
products and rejection probabilities to be evaluated using rational
arithmetic.  Only $N$ changes between the two runs.

\paragraph{Cutoffs and reconstruction.}
Both runs use four complete reflection periods and a terminal size cutoff
of twenty:
\begin{equation}\label{eq:fixed-qt-cutoffs}
K=4,\qquad M=20,\qquad z=uv=\frac1{25},\qquad
r=(uv)^2=\frac1{625}.
\end{equation}
The small boundary parameters make the edge scale decrease by a factor
of $1/625$ in each period.  By Lemma~\ref{lem:period}, four periods contain
$8N$ individual boundary reflections, in addition to the adjacent swaps.
The initial common partition is drawn from
\begin{equation}\label{eq:fixed-qt-seed}
\rho_{20}(\lambda)
=\frac{\ind_{\{\size{\lambda}\leq20\}}
        \beta_\el(\lambda)^2\,25^{-\size{\lambda}}/b_\lambda}
       {\displaystyle\sum_{\size{\mu}\leq20}
        \beta_\el(\mu)^2\,25^{-\size{\mu}}/b_\mu}.
\end{equation}
This is \eqref{eq:seed-trunc} for the parameters above.  It is not obtained
by fixing the initial partition to be empty.

The recorded history is then replayed backwards.  Bulk and boundary
increments use $G_{L,L}(xy)$ and $B_{\el,L}(hx)$, respectively, with the
\emph{recorded pre-move parameters}; reflected edge parameters need not
equal their original values in \eqref{eq:fixed-qt-word}.  The new shapes
are sampled from \eqref{eq:bulk-shape} and \eqref{eq:boundary-shape}, and
same-sign moves use the two full pre-swap edge weights.  Consequently,
the output has law $\widehat\pi_{4,20}$ in \eqref{eq:output-law} for the
chosen word $w_N$.  Exactness here refers to this finite-approximation
law; its difference from $\pi_{w_N}$ is controlled by the truncation
bounds.  The restriction $\size{\Lambda}\leq20$ applies only to the
terminal seed, not to the reconstructed partitions.  The recorded runs
initialize the pseudorandom generator with seed \texttt{20260928} at each
value of $N$; this numerical seed is distinct from the random partition
$\Lambda$.

\paragraph{Error for this parameter choice.}
The terminal estimate can be sharpened for the present $\el/\el$ example.
Indeed, writing $w=q^{a_\lambda(s)}t^{\ell_\lambda(s)}\in(0,1]$ gives
\[
\frac12=\frac{1-t}{1-q}
\leq b_\lambda(s)=\frac{1-tw}{1-qw}\leq1.
\]
Hence $\beta_\el(\lambda)\leq1$ and
$h_{\el,\el}(\lambda)\leq2^{\size{\lambda}}$.
Since $\seedZ(1/25)\geq1$, the partition-count bound
\eqref{eq:partition-count} implies
\begin{equation}\label{eq:fixed-qt-tail}
\begin{split}
\rho(\size{\Lambda}>20)
&\leq\sum_{n>20}p(n)(2/25)^n\\
&\leq
\exp\!\left(\frac{\pi^2}{6\log(25\eta/2)}\right)
\frac{\eta^{21}}{1-\eta},
\qquad \eta=\sqrt{2/25}.
\end{split}
\end{equation}
Here $p(n)$ is the number of partitions of $n$.  The last inequality uses
$\pi\sqrt{2n/3}\leq\delta n+\pi^2/(6\delta)$ with
$\delta=\log(25\eta/2)>0$, followed by a geometric sum.
The last expression in \eqref{eq:fixed-qt-tail} is approximately
$1.55782\times10^{-11}$.
For \eqref{eq:fixed-qt-parameters}--\eqref{eq:fixed-qt-word}, the constants
in \eqref{eq:explicit-bound-constants} simplify to
\[
A_*=\frac{18N}{11},\qquad B_*=\frac{45N(N-1)}8.
\]
Combining \eqref{eq:reflection-depth-bound} with
\eqref{eq:fixed-qt-tail} bounds the total variation error by
$\mathcal H((1/625)^4)$ plus the terminal term.  The evaluated bounds,
rounded upwards, are $1.05\times10^{-9}$ for $N=96$ and
$1.40\times10^{-9}$ for $N=128$.  Before rounding upwards, the latter
bound evaluates to approximately $1.39046\times10^{-9}$, shown as
$1.39\times10^{-9}$ in the figure.
These are truncation-error bounds, not errors estimated from the samples
or bounds on their distance from a limit shape.

\paragraph{Reading the figure.}
The sampled partition sequence is converted deterministically into a
rail-yard dimer configuration using the correspondence described in
\cite[Sections~2--3]{LiShi2026}.  In the ordinary geometric coordinate,
partition $\lambda^{(i)}$ has particles at
\[
y_j^{(i)}=\lambda_j^{(i)}-j+\frac12,\qquad j\geq1,
\]
where parts beyond its length are zero.  The two upper panels show
locally averaged particle occupancy for one sample at each value of $N$.
They are spatial averages of individual configurations, not averages
over repeated samples.  The plotted averaging windows contain nine
adjacent partition columns and eleven adjacent particle rows, with
edge-value padding at the ends of the displayed arrays.
The horizontal coordinate is $i/N$; the vertical coordinate is the
ordinary height $y$, without an $N$-dependent rescaling.

The lower-left panel plots $\size{\lambda^{(i)}}$ against $i/N$.
It records the number of boxes in each partition, not a height function,
a boundary curve, or algorithmic time.  The largest sizes in the two
samples are $96$ and $137$.  The lower-right panel shows the actual
$N=128$ matching, which has $549$ diagonal dimers.  The vertical crop
only restricts the display; it imposes no additional sampling cutoff.
Corollary~\ref{cor:observables} gives the same total variation bound after
this deterministic conversion.

The figure illustrates the output of the algorithm away from $q=t$.
Its qualitative annotations about smoothness or stability describe these
two displayed samples only.  Since there is only one sample at each size,
and the vertical coordinates have not been put on a common macroscopic
scale, this comparison does not establish concentration or a limit shape
at fixed $q,t$.

\section{Specialization}\label{sec:discussion}

When $q=t$, all $b_\lambda(s)$ equal one. The equal-type bulk increment and each unrestricted $L$-boundary increment become geometric, while a mixed-type bulk increment remains Bernoulli. For $c_l,c_r\in\{\el,\oa\}$ the seed becomes
\[
\rho(\lambda)=(z;z)_\infty z^{\size{\lambda}},
\]
which can be sampled through independent part multiplicities
$m_j\sim\Geom(z^j)$. Parity restrictions remain present when either boundary is $\deel$ or $\eoa$.

Even in this Schur specialization, the fixed-size kernels given here need not coincide with the original deterministic HH coupling in \cite{Li2023}. They produce the correct distribution, and they may be replaced by any known Schur local bijection satisfying the same weight-transport identities. The total variation proof is unchanged.

\appendix
\section{Finite product formulas for implementation}\label{app:branching}

This appendix specifies the one-variable coefficients appearing in every categorical draw. If $\mu\not\subseteq\lambda$, all skew coefficients are zero. Otherwise, let $R_{\lambda/\mu}$ and $C_{\lambda/\mu}$ be the sets of cells of $\lambda$ in rows and columns, respectively, that meet $\lambda\setminus\mu$. Set $b_\mu(s)=1$ for $s\notin\mu$.

If $\lambda/\mu$ is a horizontal strip, define
\begin{equation}\label{eq:branching-products}
\psi_{\lambda/\mu}(q,t)
=\prod_{s\in R_{\lambda/\mu}\setminus C_{\lambda/\mu}}
\frac{b_\mu(s;q,t)}{b_\lambda(s;q,t)},\qquad
\varphi_{\lambda/\mu}(q,t)
=\prod_{s\in C_{\lambda/\mu}}
\frac{b_\lambda(s;q,t)}{b_\mu(s;q,t)}.
\end{equation}
They are zero for a non-horizontal strip. Then
\begin{equation}\label{eq:coeff-implementation}
\begin{array}{c|cc}
 a&u_a(\lambda/\mu)&v_a(\lambda/\mu)\\ \hline
 L&\psi_{\lambda/\mu}(q,t)&\varphi_{\lambda/\mu}(q,t)\\[2pt]
 R&\varphi_{\lambda'/\mu'}(t,q)&\psi_{\lambda'/\mu'}(t,q).
\end{array}
\end{equation}
These formulas agree with the standard branching coefficients in
\cite[Chapter~VI, \S6]{Macdonald1995} and the explicit conventions in
\cite[(2.6)--(2.9)]{LiShi2026}. Empty products are one.

For a prescribed size $m$, one may enumerate all partitions of $m$, discard zero coefficients, compute the remaining finite weights, and make a categorical draw. For a same-sign update, enumerate the partitions between the fixed outer endpoints and use the full edge weights, including their powers of $x$ and $y$. The latter powers generally cannot be discarded, since the size of the intermediate partition is not fixed in that resampling.

A useful implementation check is
\[
v_a(\lambda/\mu)b_\mu=u_a(\lambda/\mu)b_\lambda.
\]
Further finite checks are the bulk coefficient identity \eqref{eq:bulk-coeff-id}, the boundary identity \eqref{eq:boundary-coeff-id}, and the conjugation formula \eqref{eq:F-dual}. Such checks do not replace the proofs, but catch parameter-order and parity mistakes.

\section{Subexponential weight bounds and admissibility}\label{app:growth}

We prove the seed estimate \eqref{eq:seed-coefficient-bound} and the finiteness assertion used in \eqref{eq:regime}. No closed formula for the terminal partition function is needed.

\begin{lemma}[A uniform box-product bound]\label{lem:box-growth}
Let $n=\size{\lambda}$ and let $\gamma(q,t)$ be defined by \eqref{eq:seed-envelope-constants}. Then
\begin{equation}\label{eq:J-bound}
J(\lambda):=\sum_{s\in\lambda}|\log b_\lambda(s;q,t)|
\leq\gamma(q,t)\sqrt n.
\end{equation}
Consequently,
\begin{equation}\label{eq:h-bound}
0\leq h_{c_l,c_r}(\lambda)\leq e^{\gamma(q,t)\sqrt n}.
\end{equation}
\end{lemma}
\begin{proof}
Put $\vartheta=\max(q,t)$. For a cell with arm $a$ and leg $\ell$, write $w=q^a t^\ell$. Then
\[
|\log b_\lambda(s)|
=|\log(1-tw)-\log(1-qw)|
\leq\frac{|q-t|}{1-\vartheta}q^a t^\ell.
\]
Let $d(\lambda)$ be the number of distinct positive row lengths. For each fixed pair $(a,\ell)$, there are at most $d(\lambda)$ cells with those arm and leg lengths. Indeed, within a group of rows of the same length, the arm fixes the column, and the leg then fixes the row. Hence
\[
\sum_{s\in\lambda}q^{a_\lambda(s)}t^{\ell_\lambda(s)}
\leq\frac{d(\lambda)}{(1-q)(1-t)}.
\]
Since $n\geq d(\lambda)(d(\lambda)+1)/2$, we have $d(\lambda)\leq\sqrt{2n}$. This proves \eqref{eq:J-bound}.

On an allowed seed partition, each cell contributes to
$\beta_{c_l}(\lambda)\beta_{c_r}(\lambda)/b_\lambda$ with exponent $-1$, $0$, or $1$. Its logarithm is therefore bounded above by $J(\lambda)$. Disallowed partitions have zero weight, giving \eqref{eq:h-bound} in all cases.
\end{proof}

For completeness, an elementary bound on the number $p(n)$ of partitions of $n$ is
\begin{equation}\label{eq:partition-count}
p(n)\leq\exp\left(\pi\sqrt{\frac{2n}{3}}\right).
\end{equation}
For $n>0$ and $\delta>0$, the multiplicity generating function gives
\[
p(n)e^{-\delta n}\leq\prod_{j\geq1}(1-e^{-\delta j})^{-1}.
\]
Taking logarithms and using $e^{\delta m}-1\geq\delta m$ yields
\[
\log p(n)\leq\delta n+
\sum_{m\geq1}\frac1{m(e^{\delta m}-1)}
\leq\delta n+\frac{\pi^2}{6\delta}.
\]
Choose $\delta=\pi/\sqrt{6n}$. The case $n=0$ is immediate.
Combining \eqref{eq:h-bound} and \eqref{eq:partition-count} proves \eqref{eq:seed-coefficient-bound}. In particular, $\seedZ(z)<\infty$ for every $0<z<1$, and \eqref{eq:seed-tail-bound} follows as stated.

\begin{proposition}[Finiteness in the subunit regime]\label{prop:finiteZ}
Under \eqref{eq:regime}, $1\leq Z_w<\infty$ for every finite labelled word and every pair of boundary types in $\cset$.
\end{proposition}
\begin{proof}
Any nonzero branching coefficient in \eqref{eq:branching-products} is a product of $b$-ratios over subsets of cells. Therefore
\[
u_a(\lambda/\mu),\ v_a(\lambda/\mu)
\leq\exp\bigl(\gamma(q,t)(\sqrt{\size{\lambda}}+\sqrt{\size{\mu}})\bigr).
\]
The same bound applies to $R$ after conjugation and parameter exchange, since $\gamma(q,t)=\gamma(t,q)$. Likewise, each nonzero $\beta_c(\lambda)$ and $\tau_c(\lambda)$ is at most $e^{\gamma(q,t)\sqrt{\size{\lambda}}}$.

Write $m_i=\size{\lambda^{(i)}}$, $d_i=|m_i-m_{i-1}|$, and $D=\sum_i d_i$. The total coefficient in \eqref{eq:target-weight}, excluding the monomials in $u,v,x_i$, is at most
\[
\exp\left(2\gamma(q,t)\sum_{i=0}^N\sqrt{m_i}\right).
\]
For fixed sizes there are at most $\prod_i p(m_i)$ partition sequences. Also $m_i\leq m_0+D$. The signs of the word and the values of $m_0,d_1,\ldots,d_N$ determine all sizes; invalid size sequences can be discarded. For fixed $D$, the number of choices of the nonnegative $d_i$ is $\binom{D+N-1}{N-1}$.

Using $v^{m_N}\leq1$ and $x_*=\max_i x_i<1$, we obtain
\begin{equation}\label{eq:Z-majorant}
Z_w\leq\sum_{m_0,D\geq0}\binom{D+N-1}{N-1}
u^{m_0}x_*^D\exp\bigl(A\sqrt{m_0+D}\bigr),
\end{equation}
where $A=(N+1)(2\gamma(q,t)+\pi\sqrt{2/3})$. Choose
$0<\delta<\min\{-\log u,-\log x_*\}$. Since
$A\sqrt n\leq\delta n+A^2/(4\delta)$, the right side of \eqref{eq:Z-majorant} is at most
\[
\frac{\exp(A^2/(4\delta))}{(1-ue^\delta)(1-x_*e^\delta)^N}<\infty.
\]
Finally, the all-empty path has weight one, so $Z_w\geq1$.
\end{proof}

\section*{Acknowledgements}

The author thanks Vadim Gorin for posing the questions answered in this
paper. This work was supported by the National Science Foundation under
grant no.~2512154 and by the Simons Foundation under grant
no.~SFI-MPS-TSM-00025603.

\paragraph{AI usage disclosure.}
The author used ChatGPT to assist with polishing the manuscript and with
implementing the algorithm in Section~\ref{subsec:fixed-qt-example}.
All original ideas in this paper are due to the author, who takes full
responsibility for the content of the paper and the algorithmic
implementation.

\end{document}